\documentclass{amsart}

\usepackage{amssymb}
\usepackage{amscd}
\usepackage[all]{xy}

\DeclareMathOperator{\Hom}{{\rm Hom}}

\DeclareMathOperator{\Fr}{\rm Fr}
\DeclareMathOperator{\Gal}{\rm Gal}

\newcommand{\D}{\varDelta}

\newcommand{\bbz}{\mathbf z}
\newcommand{\ZZ}{\mathbb Z}

\newcommand{\QQ}{{\mathbb Q}}

\newcommand{\FF}{{\mathbb F}}
\newcommand{\T}{{\mathcal T}}

\newcommand{\abold}{{\mathbf a}}
\newcommand{\bbold}{{\mathbf b}}
\newcommand{\cbold}{{\mathbf c}}

\newcommand{\AAA}{{\mathcal A}}
\newcommand{\BBB}{{\mathcal B}}
\newcommand{\DDD}{{\mathcal D}}
\newcommand{\LLL}{{\mathcal L}}
\newcommand{\NNN}{{\mathcal U}}
\newcommand{\co}{{\mathcal O}}

\newcommand{\K}{\mathbf K}

\newcommand{\Q}{\mathbf Q}

\newcommand{\WWW}{{\mathbb W}}

\newtheorem{theo}{Theorem}[section]
\newtheorem{Lemma}[theo]{Lemma}

\newtheorem{Proposition}[theo]{Proposition}
\newtheorem{Corollary}[theo]{Corollary }

\theoremstyle{remark}
\newtheorem{rem}[theo]{Remark}

\begin{document}

\title{The universal Kolyvagin recursion implies the Kolyvagin recursion}
\date{October 16, 2002}
\author{Yi Ouyang}
\address{Department of Mathematics, University of Toronto, 
100 St George St., Toronto,  Ontario M5S 3G3, Canada}
\email{youyang@math.toronto.edu}

\begin{abstract}
let $\NNN_z$ be the universal norm distribution and $M$ a fixed
power of prime $p$, by using the double complex method employed by 
Anderson, we study the universal Kolyvagin recursion occurred in the 
canonical basis 
in the cohomology group $H^0(G_z,\NNN_z/M\NNN_z)$ as given in
Ouyang~\cite{Ouyang4}. We furthermore show that
the universal Kolyvagin recursion implies the Kolyvagin recursion 
in the theory of Euler systems(Theorem 4.5.4, Rubin~\cite{Rubin2}).
One certainly hopes this could lead a new way to find new Euler systems.
\end{abstract}
\maketitle

\section{Introduction}
Let $\FF$ be a finite real abelian
extension of $\QQ$. Let $M$ be an odd positive integer. For
every squarefree positive integer $r$ the prime factors of which
are congruent to $1$ modulo $M$ and split completely  in $\FF$, the 
corresponding Kolyvagin class
 $\kappa_r\in\FF^{\times}/\FF^{\times\, M}$  satisfies a
remarkable and crucial recursion which for each prime number
$\ell$ dividing
$r$ determines the order of vanishing of
$\kappa_r$ at each place of $\FF$ above $\ell$ in terms of
$\kappa_{r/\ell}$. In the note~\cite{AO}, Anderson and Ouyang
gave the recursion a new and universal
interpretation with the help of the double complex method.
Namely, the recursion satisfied by Kolyvagin classes
is shown to be the specialization of a universal 
recursion independent of $\FF$ satisfied by the universal
Kolyvagin classes in the group cohomology of the
universal ordinary distribution. 

In the note~\cite{AO}, the question of whether  such kind of 
Kolyvagin recursion holds or not for the universal Euler systems 
is raised. The goal of this paper to answer this question. 

Let $X$ be a totally ordered set(in application, elements in $X$
are often prime ideals in a number field $K$). Let $Z$ be the set of
formal products of elements in $X$. Let $\co$ be the integer ring of 
a number field or a local field with characteristic $0$ and $\T$
a free $\co$-module.
For every element $z\in Z$, let $x^n$ be the $x$-part of $z$ for $x\in X$
and $x\mid z$. We associate $z$ with a group $G_z=\prod_{x\mid z}G_{x^n}$
and a $\T[G_z]$-module $\NNN_z$, 
a universal model satisfying certain distribution relations. This module 
$\NNN_z$ is
called the \emph{universal norm distribution} of level $z$, a generalization
of the  \emph{universal ordinary distribution} of 
Kubert~\cite{Kubert1} and the \emph{universal Euler system} of
Rubin~\cite{Rubin2}. 
A special case of the universal norm distribution was first
introduced in Ouyang~\cite{Ouyang3}. The more general case was defined 
and studied in Ouyang~\cite{Ouyang4}. Suppose $M$ is a nonzero element
in $\co$ dividing the order of all $G_{x^n}$, then 
the double complex method(see for example Anderson~\cite{Anderson2}
and Ouyang~\cite{Ouyang2}) produces a canonical basis for the
cohomology group $H^{\ast}(G_z,\NNN_z/M\NNN_z)$ as a free
$\T/M\T$-module. In particular, the canonical basis of 
$H^{0}(G_z,\NNN_z/M\NNN_z)$ can be given as 
$\{c_y: y\mid z, y\ \text{squarefree}\}$. Moreover, one can consider
$H^{0}(G_{z/x^n},\NNN_{z/x^n}/M\NNN_{z/x^n})$ as a subgroup of
$H^{0}(G_z,\NNN_z/M\NNN_z)$ with $\T/M\T$-basis 
$\{c_y: y\mid \frac{z}{x^n}, y\ \text{squarefree}\}$. 

In this paper, we first find a natural and well defined map $\D_x$ 
from $H^{0}(G_z,\NNN_z/M\NNN_z)$ to 
$H^{0}(G_{z/x^n}, \NNN_{z/x^n}/M\NNN_{z/x^n})$. We then show that 
$\D_x$ maps $c_{xy}$ to $c_y$ and $c_y$ to $0$ for $xy$ a squarefree
factor of $z$ by the double complex method.
Acting $\D_{x_i}$ consequently, we thus obtain a sequence
of elements in $H^{0}(G_z,\NNN_z/M\NNN_z)$. We thus say that the family
consisting of classes $\{c_y\}$
satisfies the \emph{universal Kolyvagin recursion}. 
In addition, we show that
the family of the original Kolyvagin classes(i.e., the cocycles which map 
to $D_r\xi_r$ as in Rubin~\cite{Rubin1}) also satisfies the 
universal Kolyvagin recursion.

Now given a number field $K$ and 
a $p$-adic representation $T$ of $G_K$ with coefficients in $\co$,
the universal norm distribution $\NNN_z$ in this situation becomes
the  universal Euler system, for which  Rubin introduced 
and studied extensively in Chapter $4$ of \cite{Rubin2}, resulting 
the proof of the celebrated Theorem 4.5.4 there. If let $\xi$ be 
a $G_K$-homomorphism from $\NNN_z$ to $H^1(F(z),T)$ for $F(z)$ a
certain extension of $K$. Let $W_M=M^{-1}T/T$, then $\xi$ induces
a map from $H^0(G_z,\NNN_z/M\NNN_z)$ to $H^1(F,W_M)$. 
Theorem 4.5.4 then states that the family of the Kolyvagin classes
satisfies an important recursive relation, which is essential to the 
effectiveness of Euler system. We call this recursion the \emph{Kolyvagin 
recursion}. We show that the universal Kolyvagin recursion implies 
the  Kolyvagin recursion in Theorem~\ref{theo:koly2}.

The author sincerely thanks Professor Greg W. Anderson for his 
continuous support and for his many ideas that lead to this paper. He also
deeply thank Professor Kumar Murty for his support. Part of this research 
is supported by the Ganita lab of University of Toronto.

\section{Review about the universal norm distribution}
\label{section:review}
We first give a brief review of the results in Ouyang~\cite{Ouyang4}.

\subsection{The universal norm distribution} \label{subsection:und}
Let $X$ be a given totally ordered set.
Let $Y$ be the set of all squarefree formal products of elements in 
$ X$, Let $Z$ be the set of all (finite or infinite) 
formal product of elements in $X$. Elements in $X$ will be denoted by
$x, x'$ and be called \emph{primes}, elements in $Y$ will be denoted by
$y, y'$ etc and elements in $Z$ will be denoted by $z, z', w$ and $\bbz$.
Moreover, $y$, $y'$  and $w$ will usually assumed to be \emph{finite}
(i.e., finite product of $x\in X$) and $\bbz$ will assume to be
\emph{infinite}.

For every $z\in Z$, the \emph{support} 
of $z$ is the unique element $\bar z\in Y$ such that 
if $x\mid z$ then $x\mid \bar z$. 
For every $z\in Z$, a stalk of $z$ is a factor $z'\mid z$ satisfying 
$\gcd (z', z/z')=1$ and is denoted by $z'\mid_s z$. 
Fix $z$, for each $y\mid \bar z$, let $z(y)$ be 
corresponding stalk of $z$ whose  support is $y$. Let $v_x(z)$ be
the integer $n$ such that $x^n=z(x)$. 
        
For every $z\in Z$, $G_z$ is an abelian group which is the
direct product of $G_{z(x)}$ for all $x\mid z$. Furthermore, for
every $z'\mid z$, $G_{z'}$ is a subgroup of $G_z$ and hence also a
quotient of $G_z$.
For any $z'\mid z$ and $g\in G_z$, Let $g_{z'}$ denote
the restriction of $g$ to $G_{z'}$. 
For each pair $x\in X$ and $z\in Z$, the Frobenius element
$\Fr_{x}$ is a given element in $G_z$ whose restriction to $G_x$ is the 
identity.

Let $\co$ be an integral domain and let $\Phi$ be its fractional field. 
Let $\T$ be a fixed torsion free finitely generated $\co$-module.
For each $x\in X$, a polynomial 
\[ p(x; t)\in \T[\, t\,] \] 
is chosen corresponding to $x$.

For every finite $z\in Z$, $\BBB_z$ is 
the free $\T$-module $\T[G_{z}]=\langle B_z\rangle_{\T}$
generated by 
\[ B_z=\{[g\ z]: g\in G_z\}. \] 
For every $z\in Z$, finite or not, let
\[ A_z=\bigcup_{\substack{z'\ \text{finite}\\ z'\mid_s z}}B_{z'} \]
and let
\[ \AAA_z=\langle A_z\rangle_{\T}=
\bigoplus_{z'\mid_s z} \BBB_{z'}=\langle[g\ z']: z'\mid_s z, z'\ 
\text{finite},\ g\in G_{z'}\rangle_{\T}. \] 
For every pair $z'\mid_s z$, the group $G_z$ acts on $\AAA_{z'}$ by 
\[ g\cdot [g''z'']= [g_{z''} \cdot g'' z''],\ z''\ \text{finite},
z''\mid_s z' \]
and by this way $\AAA_{z'}$ becomes a $\T[G_z]$-module. 

Let $\lambda_{z(x)}$
be the $\T[G_z]$-homomorphism given by 
\[ \lambda_{z(x)}: [z']\longmapsto \begin{cases} 
p(x; \Fr^{-1}_{x})[z']-N_{z(x)}[z(x)z'],\
&\text{if}\ x\nmid z',\\ 0,\ &\text{if}\ x\mid z'.\end{cases} \]
Let $\DDD_z$ be the submodule of $\AAA_z$ generated by the images
of $\lambda_{z(x)}(\AAA_{z/z(x)})$ for all $x\mid z$.
The \emph{universal norm distribution}
$\NNN_z$  is then defined to be the
quotient $\T[G_z]$-module $\AAA_z/\DDD_z$.

From Ouyang~\cite{Ouyang4}, for any $z'\mid_s z\in Z$, the map 
\[ \NNN_{z'}\rightarrow \NNN_z \]
induced by the inclusion $\AAA_{z'}\subseteq \AAA_z$ is an injective
$G_z$-homomorphism of free $\T$-modules with free cokernel
and hence the induced map 
\[ H^0(G,\NNN_{z'}/M\NNN_{z'})\rightarrow H^0(G,\NNN_z/M\NNN_z) \]
is also injective.
Thus we henceforth identify $\NNN_{z'}$
(resp.~$H^0(G,\NNN_{z'}/M\NNN_{z'})$) with a subgroup of $\NNN$ 
(resp.~$H^0(G,\NNN_z/M\NNN_z)$).
Note that we have
\[ \NNN_{z}=\bigcup_{\substack{z'\ \text{finite}\\ z'\mid_s z}}\NNN_{z'},
\;\;\; H^0(G,\NNN_z/M\NNN_z)=
\bigcup_{\substack{z'\ \text{finite}\\ z'\mid_s z}}
H^0(G,\NNN_{z'}/M\NNN_{z'}).  \]

\subsection{Anderson's resolution $\LLL_z$}
For every $z\in Z$, $\LLL_z$ is the free $\T$-module  generated by
\[ \{ [a,y]:  [a]\in A_{z/z(y)}, y\mid \bar z \} \]
We equip  $\LLL_z$ with a grading by declaring that
\[ \deg [a,y]= -\deg  y=-(\text{number of primes $x\mid y$}). \]
For any
$g\in G_z$ and $[g'z']\in A_{z/z(y)}$, We equip  $\LLL_z$ with
a $G_z$-action by the rule
\[ g[g'z',y]:= [g_{z'} g'z', y], \]
then $\LLL_z$ becomes a graded $\T[G_z]$-module. 
$\LLL_z$ is bounded above and is bounded if and only if $z$ is finite.
We equip $\LLL_z$ with a $G_z$-equivariant differential $d$ of degree
$1$ by the rule
\[ d[a,y]=
\sum_{x\mid y} \omega(x,y)(p(x;\Fr^{-1}_x)[a, y/x]-N_{z(x)}[az(x),y/x]) \]
where $\omega$ is as defined as
\[ (x,y)\longmapsto \begin{cases} (-1)^{\#\{x'\mid y: x'<x\}},\ 
&\text{if $x\mid y$};\\ 0,\ 
&\text{if $x\nmid y$.} \end{cases} \]
Then the complex $(\LLL^{\bullet}_z,d)$ is acyclic for degree
$n\neq 0$ and   $H^0(\LLL^{\bullet}_z, d)\cong\NNN_z$ induced 
by 
\[ \mathbf u: [a,y]\longmapsto \begin{cases} [a],\ &\text{if}\ y=\mathbf 1;\\
0,\ &\text{if}\ y\neq \mathbf 1. \end{cases} \]
We call the complex $(\LLL^{\bullet}_z, d)$
(or simply $\LLL^{\bullet}_z$)
\emph{Anderson's resolution} of the universal norm system $\NNN_z$. 

\subsection{The double complex $\K^{\bullet,\bullet}_{z}$} 
Assume that $G_{z(x)}$ is acyclic. Let $\sigma_{z(x)}$ be a generator
of $G_{z(x)}$. Let  $\K^{\bullet,\bullet}_{z}$ 
be the free graded $\T$-module with basis
\[ \{ [a,y,w]: y\mid \bar z, a\in A_{z/z(y)}, \bar w\mid \bar z \} \]
and with the double grading given by 
\[ \deg[a,y,w]=(-\deg y, \deg w). \]
where $\deg w=\sum_{x\mid w} v_x(w)$.
We equip $\K^{\bullet,\bullet}_z$ with a $\T[G_z]$-module structure by 
the rule
\[ g[a, y, w]=[ga, y, w], \forall\ g\in G_z \]
For every $x\mid z$,
we equip  $K^{\bullet,\bullet}_z$ with $G_z$-equivariant differentials 
$d_x$ of bidegree $(1,0)$ by the rule
\[  d_x[a,y,w]=\omega(x,y) (-1)^{\sum_{x'<x}v_{x'}(w)}
\left (p(x;\Fr^{-1}_x) [a,y/x,w]-N_{z(x)}[az(x),y/x,w] \right ), \]
and with $G_z$-equivariant differentials 
$\delta_x$ of bidegree $(0,1)$ by the rule
\[ \delta_x [a, y, w]= (-1)^{\sum_{x'\leq x}v_{x'}(y)}
(-1)^{\sum_{x'<x}v_{x'}(w)} a_{z(x)}[a, y, wx] \]
where $a_{z(x)}$ is equal to $1-\sigma_{z(x)}$ if $v_{x}(w)$ even
and $N_{z(x)}$ if $v_{x}(w)$ odd. Any two distinct differentials in
the family $\{d_x\}\cup\{\delta_x\}$ anticommute. If let
\[ d=\sum_x d_x,\qquad \delta=\sum_x \delta_x, \]
then $d$ is a differential of bidegree $(1,0)$ and $\delta$ is of
bidegree $(0,1)$ and $d$ and $\delta$ are anticommute. $(\K_z; d,\delta)$
is a double complex of $G_z$-modules. Let $\K^{\bullet}_{z}$ be 
the single total complex of this double complex.

Let $\bar \K^{\bullet}_{z}$ be
the quotient of free $\T$-module generated by
\[ \{[a,w], a\in A_z,\bar w\mid \bar z \} \]
modulo relations generated by
\[ p(x;\Fr^{-1}_x)[a,w]-N_{z(x)}[a\; z(x),w],
\ a\in A_{z/z(x)},\ \bar w\mid \bar z,\  x\mid z. \]
With the grading given by the rule
\[ \deg [a,w]=\deg w \]
and the differential $\delta$ given by the rule
\[ \delta [a,w]=\sum_{x\in z} (-1)^{\sum_{x'<x}v_{x'}w}
a_{z(x)} [wx], \]
$\bar \K^{\bullet}_{z}$ is a complex of $G_z$-modules whose
cohomology is nothing but the group $H^{\ast}(G_z,\NNN_z)$.
The homomorphism
\[ \mathbf u: \K^{\bullet}_z\longrightarrow
\bar \K^{\bullet}_z,\  [a,y,w]\longmapsto
\begin{cases} [a,w],\ &\text{if}\ y=\mathbf 1 \\
0,\ &\text{if}\ y\neq \mathbf 1 \end{cases} \]
is a quasi-isomorphism. Thus it induces isomorphisms between
$H^{\ast}(\K^{\bullet}_{z}, d+\delta)$ and $H^{\ast}(G_z, \NNN_z)$(resp.
$H^{\ast}( \K^{\bullet}_{z}/M\K^{\bullet}_{z})$ and 
$H^{\ast}(G_z, \NNN_z/M\NNN_z)$ for $0\neq M\in \co$). In particular,
for any $0$-cocycle $c$ in $\K^{\bullet}_{z}/M\K^{\bullet}_{z}$, the map
$\mathbf u$ sends its bidegree $(0,0)$-component 
$\sum [a,\mathbf 1,\mathbf 1]$ to $\sum [a,1]\in 
\bar \K^{0}_{z}/M\bar \K^{0}_{z}$ and then to $\sum [a]\in \NNN_z/M\NNN_z$,
the resulting element is fixed by $G_z$ and hence is a cocycle
in $H^0(G_z,\NNN_z/\NNN_z)$.

For every $z'\mid_s z$, let $\K_{z'}$ be the 
submodule of $\K_{z}$ generated by 
\[ \{[a, y, w]: y\mid z', a\in B_{z'/z(y)}, 
\bar w\mid z'\} \]
and let $\K_{z}(z')$ be the submodule generated by
\[ \{[a,y,w]:y\mid z', a\in B_{z'/z(y)}, 
\bar w\mid z\} \]
Then $\K_{z'}$ and $\K_{z}(z')$ are compatible
with differentials $d_x$ and $\delta_x$. The $(d+\delta)$-cohomology 
of $\K_{z'}$(resp. $\K_{z'}/M\K_{z'}$) is 
just $H^{\ast}(G_{z'}, \NNN_{z'})$(resp. 
$H^{\ast}(G_{z'}, \NNN_{z'}/M\NNN_{z'})$) 
and the $(d+\delta)$-cohomology
of $\K_{z}(z')$(resp. $\K_{z}(z')/M\K_{z}(z')$)
 is $H^{\ast}(G_z, \NNN_{z'})$(resp. 
$H^{\ast}(G_z, \NNN_{z'}/M \NNN_{z'})$).

\subsection{The canonical basis for $H^{0}(G_z,\NNN_z/M\NNN_z)$} 
Suppose now that $M$ is a common divisor of $|G_{z(x)}|$ 
and $p(x;1)$ for every $x\mid z$.
Let $\mathbf S_z$ be the $\T$-submodule of $\K_{z}$ generated by 
\[ \{[a, y, w]: a\in B_{z/z(y)},\ y\mid z,\ 
\bar w\mid z, a\notin B_{\mathbf 1}\ \text{if}\ y\mid w \}.
\]
Then
$\mathbf S_z/M\mathbf S_z$ is a submodule of $\K_{z}/M\K_z$ 
which is also $d-$ and $\delta-$ stable and thus is a subcomplex of 
$\K_{z}/M\K_z$ with respect to the multi-complex structure of $\K_z/M\K_z$. 
Let $\Q_{z}/M\Q_z$ be the quotient of $\K_z/M\K_z$ by 
$\mathbf S_z/M\mathbf S_z$.
Then $\Q_{z}/M\Q_z$ is a free $\T/M\T$-module generated by
\[ \{[\mathbf 1, y, w]: y\mid \bar w\mid z \} \] 
with all induced differentials $d=\delta=0$.  Write 
the quotient map from $\K_z/M\K_z$ to $\Q_z/M\Q_z$ as $\rho_M$.
The homomorphism $\rho_M$ is a quasi-isomorphism as claimed in 
Ouyang~\cite{Ouyang4}. Every element $[\mathbf 1, y,y]$ for
$y$ finite and dividing $z$ in $\Q_z/M\Q_z$ thus induces a 
$0$-cocycle in $\K_z/M\K_z$ and henceforth induces a cocycle
in $H^0(G_z,\NNN_z/M\NNN_z)$. Write this cocycle $\bar c_y$, then
\[ \{ \bar c_y: y\ \text{finite}, y\mid z\} \]
forms a basis for $H^0(G_z,\NNN_z/M\NNN_z)$. We call it
the \emph{canonical basis}.

\section{The universal Kolyvagin Recursions} \label{section:recursion}
\subsection{The Kolyvagin conditions}\label{subsection:KolyvaginC}
In this section, we fix an infinite $\bbz\in Z$. we write
$\NNN_{\bbz}$($\LLL_{\bbz}$, $G_{\bbz}$ and $\K_{\bbz}$ etc.) 
as $\NNN$($\LLL$, $G$ and $\K$ etc.). For any $x\mid \bbz$, let 
$r_x(t)\in \co[t]$ and 
\[ \gamma_{\bbz(x)}=p(x;\Fr^{-1}_x)-r_x(\Fr^{-1}_x)|G_{\bbz(x)}|. \]
Assume the 
following \emph{Kolyvagin conditions} hold:
\begin{itemize}
\item There exists $M\in \co$ such that $M\mid |G_x|$ and $M\mid p(x;1)$
for every $x\mid \bbz$. We fix $M$ here after.
\item The group $G_{\bbz(x)}$ is cyclic for every $x\mid \bbz$.
\item The homomorphism $\gamma_{\bbz(x)}: \NNN_{\bbz/z(x)}\rightarrow 
\NNN_{\bbz/z(x)}$ 
has a trivial kernel for every $x\mid \bbz$.
\end{itemize}
Following the first two assumptions, we can and hence will apply the results
in \S~\ref{section:review}. If we let
$n_{x}=\textrm{Norm}(x)$ and let 
\[ r_x(\Fr^{-1}_x)=\frac{p(x;\Fr^{-1}_x)-p(x;n_x\Fr^{-1}_x)}{|G_{\bbz(x)}|}, \]
then $\gamma_{\bbz(x)}=p(x;n_x\Fr^{-1}_x)$. When $p(x;t)$ is coming from
certain characteristic polynomial  in a $p$-adic representation,
then  $\gamma_{\bbz(x)}$ can be shown to satisfy the last assumption.
See Lemma~\ref{lemma:gammax} in \S~\ref{section:kolyeuler} for more details.

\subsection{The submodule $I_x$ of $\NNN$} Let $x\mid \bbz$ be given.
 We define
\[ I_x\subset \NNN \]
to be the $\T[G]$-submodule generated by all elements of $\NNN$ represented by 
the expressions of the form
\[ [z\bbz(x)]-g[z\bbz(x)]\ \ \text{or}\ \  
r_x(\Fr^{-1}_x) [z]-g[z\bbz(x)]
(g\in G_{\bbz(x)}, z\ \text{finite}\ z\mid_s \bbz, x\nmid z ). \]
When the dependence on $\bbz$ is needed to emphasize, 
we write $I_x$ as $I_{x,\bbz}$.
Note that 
\[ (\sigma_{\bbz(x)}-1)\NNN \subset I_x, \] 
the quotient $\NNN/I_x$ can be viewed as a $G/G_{\bbz(x)}$-module.

\begin{Proposition} \label{Proposition:exact}
The sequence
\[ 0\longrightarrow \NNN_{\bbz/\bbz(x)}
\xrightarrow{\gamma_{\bbz(x)}}
\NNN_{\bbz/\bbz(x)}
\longrightarrow \NNN/I_x
\longrightarrow 0 \]
is exact where the map $\NNN_{\bbz/\bbz(x)}\rightarrow \NNN/I_x$
is that induced by the inclusion of $\NNN_{\bbz/\bbz(x)}\subset \NNN$.
\end{Proposition}
\begin{proof}
We consider the complex homomorphism
\[ \gamma_{\bbz(x)}: \LLL^{\bullet}_{\bbz/\bbz(x)}\longrightarrow 
\LLL^{\bullet}_{\bbz/\bbz(x)} \]
The mapping cone of  $\gamma_{\bbz(x)}$ is just the complex
\[ Cone^{\bullet}(\gamma_{\bbz(x)})=
(\LLL^{n+1}_{\bbz/\bbz(x)},\LLL^{n}_{\bbz/\bbz(x)}) \]
with the differential given by
\[ d(a,b)=(-da, \gamma_{\bbz(x)}(a)+db), \text{for}\ (a,b)
\in Cone^{n}(\gamma_{\bbz(x)}). \]
From homological algebra, we know there exists 
the following exact sequence
\[ \Sigma: 0\longrightarrow \LLL^{\bullet}_{\bbz/\bbz(x)}\longrightarrow 
Cone^{\bullet}(\gamma_{\bbz(x)})\longrightarrow  
\LLL^{\bullet}_{\bbz/\bbz(x)}[1]
\longrightarrow 0. \]
Let
\[ s_x: \LLL\rightarrow \LLL(\bbz/\bbz(x)) \]
be the unique homomorphism such that
\[ s_x[a, y]\equiv \begin{cases}
\omega(x,y)[a, y/x]\ &\mbox{if
$x\mid  y$}\\ 0\ &\mbox{otherwise} \end{cases}  \]
for all symbols $[a,y]$ in the canonical basis of $\LLL$.
The homomorphism $s_x$ is of degree $1$ and satisfies the relation
\[ s_x d=-ds_x \]
as can be verified by a straightforward calculation. Now consider the
sequence
\[ \Sigma':\;0\rightarrow \LLL^{\bullet}_{\bbz/\bbz(x)}\longrightarrow
\LLL^{\bullet}/\LLL'\stackrel{s_x}\rightarrow 
\LLL^{\bullet}_{\bbz/\bbz(x)}[1]\rightarrow 0 \]
where $\LLL'$ is the $\T[G]$-submodule of $\LLL$
generated by all elements of the form
\[ [z\bbz(x), y]-[gz\bbz(x), y]\ \ \text{or}\ \  
r_x(\Fr^{-1}_x)[z, y]-[gz\bbz(x), y]
( \bbz(y)z\bbz(x)\mid \bbz,\ g\in G_{\bbz(x)}). \]
and the map $\LLL_{\bbz/\bbz(x)}\rightarrow \LLL/\LLL'$ is that induced
by the inclusion $\LLL_{\bbz/\bbz(x)}\subset \LLL$. It is easy to verify that 
$\Sigma'$ is short exact. The two complexes $\Sigma$ and $\Sigma'$ are 
actually isomorphic: just let the two side maps be the identities and 
let the middle map from $\LLL^{\bullet}/\LLL'$ to $Cone^{\bullet}$ be given by
\[ [a, yx]\mapsto (\omega(x,xy)[a,y],0),
\qquad [a, y]\mapsto (0, [a,y]), \forall\ y\mid \bbz/\bbz(x). \]
Since $\LLL'$ is a graded $d$- and $G$-stable subgroup of $\LLL^{\bullet}$
and $(\sigma_{\bbz(x)}-1)\LLL\subset \LLL'$,
it follows that $\Sigma'$ can be viewed as a
short exact sequence of complexes of
$G/G_{\bbz}$-modules.  Because 
$H^{\ast}(\LLL^{\bullet}_{\bbz/\bbz(x)},d)$ is 
concentrated in degree
$0$, the long exact sequence of
$G/G_{\bbz(x)}$-modules deduced from
$\Sigma'$ by taking $d$-cohomology has at most four nonzero terms
and after making the evident identifications takes the form
\[ \dots\rightarrow 0\rightarrow 
H^{-1}(\LLL^{\bullet}/\LLL',d)\rightarrow
\NNN_{\bbz/\bbz(x)}\xrightarrow{\gamma_{\bbz(x)}}
\NNN_{\bbz/\bbz(x)}\longrightarrow
\NNN/I_x \rightarrow 0\rightarrow \dots \]
where the map $\NNN_{\bbz/\bbz(x)}\rightarrow \NNN/I_x$ is that induced by the
inclusion $\NNN_{\bbz/\bbz(x)}\subset \NNN$. By the assumption of
the Kolyvagin conditions, we have 
\[ H^{-1}(\LLL/\LLL',d)=\ker\left(\NNN_{\bbz/\bbz(x)}
\xrightarrow{\gamma_{\bbz(x)}} \NNN_{\bbz/\bbz(x)}\right)=0, \]
whence the result.
\end{proof}

\begin{theo}\label{PropDef:DeltaDef}
For every prime number $x$ dividing $\bbz$ there exists a unique
homomorphism
\[ \D_x: H^0(G,\NNN/M\NNN)\rightarrow
H^0(G,\NNN_{\bbz/\bbz(x)}/M\NNN_{\bbz/\bbz(x)}) \]
such that
\[ \frac{(1-\sigma_{\bbz(x)} )a}{M}\equiv
\frac{\gamma_x b}{M}\bmod{I_x}\Leftrightarrow
 \D_x(a\bmod{M\NNN})=b \bmod{M\NNN_{\bbz/\bbz(x)}} \]
for all
$a\in \NNN$ representing a class in $H^0(G,\NNN/M\NNN)$
and
$b\in \NNN_{\bbz/\bbz(x)}$ representing a class in
$H^0(G,\NNN_{\bbz/\bbz(x)}/M\NNN_{\bbz/\bbz(x)})$.  Moreover one has
\[ \D_x H^0(G,\NNN_z/M\NNN_z)\subset 
H^0(G,\NNN_{z/\bbz(x)}/M\NNN_{z/\bbz(x)}) \]
for all finite $z\mid_s \bbz$ and divisible by $x$.
\end{theo}
\begin{proof}  Put
\[
\begin{array}{rcl}
A&:=&\left\{a\in \NNN\mid \mbox{$a$ represents a class in
$H^0(G,\NNN/M\NNN)$}\right\},\\
B&:=&\left\{b\in \NNN_{\bbz/\bbz(x)}\mid \gamma_x(b)\in
M\NNN_{\bbz/\bbz(x)}\right\},\\
C&:=&\left\{(a,b)\in X\times Y\left|
\frac{(1-\sigma_{\bbz(x)})a}{M}\equiv\frac{\gamma_x(b)}{M}\bmod{I_x}
\right.\right\}.
\end{array}
\]
Fix a finite $z\mid_s \bbz$ and divisible by $x$.
To prove the proposition it is enough to prove the following three
claims:
\begin{enumerate}
\item $C\cap(M\NNN\times B)=M\NNN\times M\NNN_{\bbz/\bbz(x)}$.
\item $(\sigma-1)C\subset M\NNN\times M\NNN_{\bbz/\bbz(x)}$ 
for all $\sigma\in G$.
\item For all $a\in A\cap \NNN_{z}$ there exists  
$b\in B\cap \NNN_{z/\bbz(x)}$ such that $(a,b)\in C$.
\end{enumerate}

We turn to the proof of the first claim. 
Only the containment
$\subset$ requires proof; 
the containment $\supset$ is trivial. Suppose
we are given 
$(a,b)\in C\cap(M\NNN\times B)$. Then
$\frac{\gamma_x(b)}{M}\in I_x\cap
\NNN_{\bbz/\bbz(x)}$ and hence by
Proposition~\ref{Proposition:exact} there exists
$c\in \NNN_{\bbz/\bbz(x)}$ such that
$\gamma_x(b)=M\gamma_x(c)$.
It follows that
$b=Mc$. Thus the first claim is proved. The second
claim follows immediately from the first.

We turn finally to the proof of the third
claim. Let
\[ \beta_{\bbz(x)}:\AAA_{z}\rightarrow\AAA_{z/\bbz(x)} \]
be the unique homomorphism such that
\[ \beta_{\bbz(x)}[z']=r_x(\Fr^{-1}_x)[z'/\bbz(x) ] \]
for all $x\mid z'\mid_s z$. 
For each prime number $x$ dividing $z$, recall that
\[ \lambda_{\bbz(x)}:\AAA_{z/\bbz(x)}\rightarrow\AAA_z \]
is the unique homomorphism such that
\[ \lambda_{\bbz(x)}[z']:=p(x;\Fr^{-1}_x)[z']- N_{\bbz(x)}[z'\bbz(x)] \]
for all $x\nmid z'$. 
Note that $\beta_{\bbz(x)}$ commutes with $\lambda_{\bbz(x')}$ for 
$x'\neq x$ and that 
the composite homomorphism $\beta_{\bbz(x)}\lambda_{\bbz(x)}$ induces the
endomorphism $\gamma_{\bbz(x)}$ of $\AAA_{z/\bbz(x)}$.
Choose a lifting
$\mathbf a\in \AAA_z$ of $a$. By hypothesis there exists an identity
\[ (\sigma_{\bbz(x)}-1)\mathbf a=M\mathbf b+\sum_{x\mid z}\lambda_{\bbz(x)}
\mathbf b_x\;\;\;(\mathbf b\in \AAA_z,\;\; \mathbf b_x\in \AAA_{z/\bbz(x)}), \]
and hence also an identity
\[ 0=M\beta_{\bbz(x)}\mathbf b
+\gamma_{\bbz(x)}\mathbf b_x+
\sum_{x'\mid \frac{z}{\bbz(x)}}\lambda_{\bbz(x')}
\beta_{\bbz(x)}\mathbf b_{x'}. \]
Then the element
$b\in \NNN_{z/\bbz(x)}$ represented by
$\mathbf b_x$ has the desired property, namely that
$(a,b)\in C$. Thus the third claim is proved and with it the
result.
\end{proof}

\subsection{The universal Kolyvagin recursion}
We say that a family of classes
\[ \{c_y\in H^0(G,\NNN/M\NNN)\}_{y\mid \bar \bbz} \]
indexed by finite $y\mid \bar \bbz$  satisfies the {\em
universal Kolyvagin recursion} if the following
conditions hold for all finite $y\mid \bar \bbz$ and primes 
$x\mid \bbz$:
\begin{itemize}
\item $c_y\in
H^0(G_{\bbz(y)},\NNN_{\bbz(y)}/M\NNN_{\bbz(y)})=
H^0(G,\NNN_{\bbz(y)}/M\NNN_{\bbz(y)})\subset
H^0(G,\NNN/M\NNN)$.
\item $x\mid y\Rightarrow  \Delta_x c_y=c_{y/x}$.
\end{itemize}

\subsection{The diagonal shift operator $\Delta_x$} 
For each  $x$ dividing
$\bbz$, we define the corresponding {\em diagonal
shift} operator
$\Delta_x$ on $\K$ of bidegree $(1,-1)$   by the
rule
\[ \Delta_x[a,y,w]:=
\begin{cases}
[a,y/x,w/x],\  &\text{if $x\mid y$ and $x\mid w$,}\\
0,\ &\text{otherwise.}
\end{cases}
\]
One has
\[ \Delta_x d_{x'}=d_{x'}\Delta_x,\;\;\;
\Delta_{x}\delta_{x'}=\delta_{x'}\Delta_x \]
for all primes $x'\mid \bbz$ distinct from $x$.
One has
\[ \Delta_x d_x=d_x \Delta_x=0,\;\;\;
(\delta_x
\Delta_x-\Delta_x\delta_x)\K\subset M\K. \]
For every finite $z$ dividing $\bbz$ one has
\[ \Delta_x \K_z\subset \begin{cases}
\K_{z/\bbz(x)},\ &\text{if $x\mid z$,}\\
\{0\},\ &\text{otherwise.}
\end{cases} \]
The action of
$\Delta_x$ therefore passes to 
\[ H^0(\K_z/M\K_z,d+\delta)=H^0(G_z,\NNN_z/M\NNN_z) \]
and in the limit to
\[ H^0(\K/M\K,d+\delta)=H^0(G,\NNN/M\NNN). \]

\begin{Proposition}
For each $x\mid \bbz$, the
endomorphism of 
$H^0(G,\NNN/M\NNN)$ induced by the diagonal shift
operation $\Delta_x$ coincides with
$\Delta_x$ defined in Theorem~\ref{PropDef:DeltaDef}. 
\end{Proposition}
\begin{proof} Fix a finite $z\mid_s \bbz$ divisible by $x$.
Fix a class
\[ c\in H^0(G_z,\NNN_z/M\NNN_z).\]
It suffices to show that $\D_x$ and 
the endomorphism of $H^0(G,U/MU)$ induced by $\Delta_x$
applied to $c$ give the same result.
Let
$\cbold$ be a
$0$-chain in $\K_z$ reducing modulo $M\K_z$ to a
$0$-cycle representing $c$. Write
\[ 0=(d+\delta)\cbold+M\bbold \]
where $\bbold$ is a $1$-chain of $\K_z$.
For any finite $y\mid z$
and finite $w$ such that $\bar w\mid z$, let 
\[ (\abold\mapsto \abold\otimes [y,w]):\AAA_{z/\bbz(y)}
\rightarrow \K_z \]
be the unique homomorphism such that
\[ [a]\otimes [y,w]:=[a,y,w] \]
for all $a\in A_{z/\bbz(y)}$.  Write
\[ \cbold=\sum \cbold_{y,w}\otimes
[y,w],\;\;\;\Delta_x\cbold=\sum
\cbold_{yx,wx}\otimes
[y,w]\;\;\;(\cbold_{y,w}\in
\AAA_{z/\bbz(y)})
\]
and
\[ \bbold=\sum\bbold_{y,w}\otimes[y,w]\;\;\;
(\bbold_{y,w}\in \AAA_{z/\bbz(y)}) \]
where all the sums are extended over pairs $(y,w)$
consisting of finite $y\mid z$
and $w$ with $\bar w\mid z$. Let
$\beta_{\bbz(x)}:\AAA_z\rightarrow\AAA_{z/\bbz(x)}$
be as in the proof of
Proposition~\ref{PropDef:DeltaDef}.
By hypothesis one has an identity
\[ 0=\sum_{\substack{x'\mid z\\ x'<x}}\lambda_{\bbz(x')}\cbold_{x',x}
+\lambda_{\bbz(x)}\cbold_{x,x}
-\sum_{\substack{x'\mid z\\ x'>x}} \lambda_{\bbz(x')}\cbold_{x',x}
+(1-\sigma_{\bbz(x)})\cbold_{\mathbf 1,\mathbf 1}
+M\bbold_{\mathbf 1,x} \]
and hence also an identity
\[ 0=
\sum_{\substack{x'\mid z\\ x'<x}}
\lambda_{\bbz(x')}\beta_{\bbz(x)}\cbold_{x',x}
+\gamma_{\bbz(x)} \cbold_{x,x}
-\sum_{\substack{x'\mid z\\ x'>x}}
\lambda_{\bbz(x')}\beta_{\bbz(x)}\cbold_{x',x}
+M\beta_{\bbz(x)}\bbold_{\mathbf 1,x}. \]
Let $a\in \NNN_r$ be the element represented by
$\cbold_{\mathbf 1,\mathbf 1}$ and let $b\in \NNN_{z/\bbz(x)}$ be the
element represented by
$\cbold_{x,x}$. 
One the one hand, the class of $H^0(G_z,\NNN_z/M\NNN_z)$ represented by 
the $0$-cocycle $\cbold \bmod{M}$ of the complex $(\K_z/M\K_z,d+\delta)$
is $a\bmod{M\NNN_z}$ and the class of
$H^0(G_{z/\bbz(x)},\NNN_{z/\bbz(x)}/M\NNN_{z/\bbz(x)})$ represented by the
$0$-cocycle $\Delta_x \cbold\bmod{M}$ of the complex
$(K_{z/\bbz(x)}/MK_{z/\bbz(x)},d+\delta)$ is $b\bmod{M\NNN_{z/\bbz(x)}}$.
But on the other hand, one has
\[ \frac{(1-\sigma_{\bbz(x)})a}{M}\equiv
\frac{\gamma_{\bbz(x)} b}{M}\bmod{I_x} \]
and hence
\[ \D_x(a\bmod M\NNN)\equiv b\bmod{M\NNN_{\bbz/\bbz(x)}} \]
by Theorem~\ref{PropDef:DeltaDef}. 
Therefore the results of applying $\D_x$ and the endomorphism of
$H^0(G,\NNN/M\NNN)$  induced by $\Delta_x$ to the class
$a\bmod{M\NNN}$ indeed coincide.
\end{proof}

\begin{Corollary}\label{Corollary:Crucial}
The canonical basis $\{\bar{c}_y: y\ \text{finite}, y\mid \bbz\}$ satisfies 
the universal Kolyvagin recursion.
\end{Corollary}
\begin{proof} Clear.
\end{proof}
\begin{Corollary}\label{Corollary:TheoremB}
Any system of classes $\{b_y\}$ satisfying the universal Kolyvagin
recursion and the normalization
$b_{\mathbf 1}=\bar{c}_{\mathbf 1}$ is a $\T/M\T$-basis of
$H^0(G,\NNN/M\NNN)$.
\end{Corollary}
\begin{proof} Fix a finite $y\mid \bbz$, let $z=\bbz(y)$. Let 
\[ y=x_1\cdots x_n \] 
be the prime factorization of $y$. One then has
\[ \D_{x_1}\cdots \D_{x_n}b_y=b_{\mathbf 1}=\bar{c}_{\mathbf 1}
=\D_{x_1}\cdots \D_{x_n}\bar{c}_y \]
and hence
\[ b_y-\bar{c}_y\in\ker\left(H^0(G_{z},\NNN_{z}/M\NNN_z)
\xrightarrow{\D_{x_1}\cdots \D_{x_n}}
H^0(G_z,\NNN_z/M\NNN_z)\right)=\bigoplus_{\substack{
y'\mid y\\ y'\neq y}}\T/M\T\cdot \bar{c}_{y'}, \]
whence the result.
\end{proof}
\subsection{The original Kolyvagin classes}
For any $x\mid \bbz$, we let
\[ D_{\bbz(x)}=\sum_{k=0}^{|G_{\bbz(x)}|-1} k\cdot \sigma^k_{\bbz(x)}. \]
Then one has
\[ (1-\sigma_{\bbz(x)})D_{\bbz(x)}= N_{\bbz(x)}-|G_{\bbz(x)}|. \]
For any finite $z\mid_s \bbz$, let
\[ D_y=\prod_{x\mid y} D_{\bbz(x)}. \]
In particular, $D_{\mathbf 1}=1$.
Let $c'_y=D_y[\bbz(y)]\in \NNN$. 
One sees that for $x\mid y$,
\[ \begin{split} (1-\sigma_{\bbz(x)})c'_y=&
(N_{\bbz(x)}-|G_{\bbz(x)}|)D_{y/x}[\bbz(y)] \\
=& p(x; \Fr^{-1}_x) D_{y/x}[\bbz(y/x)]-|G_{\bbz(x)}|D_{y/x}[\bbz(y)]\\
=& \gamma_{\bbz(x)} c'_{y/x}+|G_{\bbz(x)}| D_{y/x}
\left (r_x(\Fr^{-1}_x) [\bbz(y/x)]-[\bbz(y)]\right ). \end{split} \]
The above identity tells us two things. First, by induction,
we see  that $(1-\sigma_{\bbz(x)})c'_y\subset M\NNN_{\bbz(y)}$ for
every $x\mid y$, thus the image $\bar c'_y$ of $c'_y$ in 
$\NNN_{\bbz(y)}/M\NNN_{\bbz(y)}$ 
is fixed by $G_{\bbz(y)}$ and hence is a $0$-cocycle.
Secondly, we see that
\[ \frac{(1-\sigma_{\bbz(x)})c'_y}{M}=
\frac{\gamma_{\bbz(x)} c'_{y/x}}{M} \pmod {I_x}, \]
thus 
\[ \D_{x} \bar c'_y=\bar c'_{y/x}. \]
Hence $\{\bar c'_y: y\ \text{finite}, y\mid \bbz\}$ satisfies 
the universal Kolyvagin recursion.
In particular, one sees that
$\bar c'_{\mathbf 1}=\bar c_{\mathbf 1}$. By 
Corollary~\ref{Corollary:TheoremB}, we thus have
\begin{theo} \label{theo:TheoremB}
The set of classes 
\[ \{\bar c'_y=D_y [\bbz(y)] \mod M\NNN: y\ \text{finite}, y\mid \bbz\} \]
constitutes a $\T/M\T$-basis of $H^0(G,\NNN/M\NNN)$.
\end{theo}
\begin{rem} 1. The classes $\bar c'_y$ are the original classes used in the
study of Euler system. The above Theorem~\ref{theo:TheoremB} is a
generalization of Theorem B in Ouyang~\cite{Ouyang2}. The proof here
is following the proof of Theorem B given at the
end of Anderson-Ouyang~\cite{AO}.

2. Apparently the definitions of $D_x$, $D_y$ and $c'_y$ depend on 
the choice of $\bbz$. We shall use $D_{\bbz(x)}$, $D_{\bbz(y)}$ and 
$c'_{\bbz(y)}$ when emphasis of the dependence is needed. 
\end{rem}

\section{The Kolyvagin recursion in Euler systems}
\label{section:kolyeuler}
In this section, we apply the results in the previous section
to show that a family satisfying  the universal Kolyvagin recursion
maps to a family satisfying the Kolyvagin recursion in an Euler system. 
We shall follow heavily
Chapter $4$ of Rubi's book~\cite{Rubin2}, which is actually the main 
motivation for this paper. Many results there will be introduced 
here without proof.  We first give a brief review of the definition of 
the universal Euler system according to Rubin.
\subsection{The universal Euler system and the Euler system}
Let $K$ be  a fixed number field. Let $p$ be a fixed rational prime number.
Let  $\Phi$ be a finite extension of $\QQ_p$ and let $\co$ be the
ring of integer of $\Phi$. Let $T$ be a $p$-adic representation of
$G_K$ with coefficients in $\co$. Assume that $T$ is unramified outside a 
finite set of primes of $K$. Let $W=(T\otimes_{\co}\Phi)/T$. Fix a nonzero
$M\in \co$ and let $W_M$ be the $M$-torsion in $W$. Let 
$\WWW_M=Ind^{G_K}_{\{1\}}W_M$. 
The exact sequence
\[ 0\rightarrow W_M\rightarrow \WWW_M \rightarrow \WWW_M/W_M \rightarrow 0\]
thus induce a canonical (surjective) map 
\[ \delta_L: ( \WWW_M/W_M)^{G_L}\rightarrow H^1(G_L, W_M) \]
for every finite extension $L$ of $K$.

Fix an ideal $\mathfrak N$ of $K$ divisible by $p$ and by all primes where 
$T$ is ramified. Let $X$ be the set of all primes $x$ of $K$ which is
prime to $\mathfrak N$ and $K(x)\neq  K({\mathbf 1})$, where $K(x)$ 
is the ray class field of $K$ modulo $x$ and $K(1)$ is the Hilbert class 
field of $K$. Then $Y$ is the set of squarefree 
products of primes in $X$ and $Z$ is the set of formal products
of primes in $X$. Let $K(x^n)$ be the ray class field of $K$ modulo 
$x^n$. Class field theory tells us that 
$G_{x^n}=\Gal(K(x^n)/K(\mathbf 1))$ is a cyclic group and
$K(x_1^{n_1})\bigcap K(x_1^{n_1})=K(\mathbf 1)$ for $x_1\neq x_2$.
Let $\sigma_{x^n}$ be a generator of $G_{x^n}$.
For every finite $z=x^{n_1}_1\cdots x_{k}^{n_k}\in Z$, let $K(z)$ be the 
composite
\[ K(z)=K(x_1^{n_1})\cdots K(x_k^{n_k}). \]
Fix a $\ZZ^d_p$-extension $K_{\infty}/K$ which no finite prime
splits completely. We write $K\subset_f F\subset K_{\infty}$ to indicate
$F/K$ a finite subextension of $K_{\infty}/K$. For 
$K\subset_f F\subset K_{\infty}$, we let
$F(z)=FK(z)$. Let 
$G_{z}=\Gal(F(z)/F(\mathbf 1))\cong \Gal(K(z)/K(\mathbf 1))$.
We see that for any $z'\mid_s z$, $G_{z}=G_{z'}\times G_{z/z'}$.

Let $\Fr_{x}$ denote a Frobenius of $x$ in $G_K$, and let
\[ p(x;t)=\det(1-\Fr_{x}^{-1} t| T^{\ast}) \in \co[t]. \]

Let $\T=\T(F)=\co[\Gal(F(\mathbf 1)/K)]$. With the above $X$, $Y$, 
$Z$, $\co$, $\Phi$
and $p(x;t)$, the corresponding universal norm distribution 
$\NNN_{F,z}$ 
(related to $F$) is called the \emph{universal Euler system} of 
level $(F,z)$. In application, we only need a certain 
infinite subset $X_{F,M}$ of $X$ with elements $x$ satisfying the
following conditions:
\begin{itemize}
\item $M\mid |G_x|$ (hence $M\mid |G_{x^n}|$);
\item $M\mid p(x;1)$;
\item $x$ splits completely in $F(\mathbf 1)/K$.
\end{itemize}
Hereafter, the pair $(F,z)$ will always mean that $z$ is finite and
every $x\mid z$ is inside $X_{F.M}$.
An \emph{Euler system} is now a $G_K$-homomorphism
\[ \xi: \bigcup_{F, z} \NNN_{F,z}\longrightarrow
\varinjlim_{F,z} H^1(F(z),T). \]
The following Proposition is crucial to the definition of Kolyvagin
classes:
\begin{Proposition} \label{proposition:df}
Suppose $\xi$ is an Euler system. Then there exists
a family of $\co[G_K]$-module maps $\{d_F\}$ such that the following 
diagrams are commutative:
\[ \begin{CD}
\NNN_{F,z}    @>{d_F}>> (\WWW_M/W_M)^{G_{F(z)}} \\
@V{\xi}VV  @V{\delta_{F(z)}}VV  \\
H^1(F(z),T) @>>> H^1(F(z), W_M)
\end{CD}\qquad
\begin{CD}
\NNN_{F',z}    @>{d_{F'}}>> (\WWW_M/W_M)^{G_{F'(z)}} \\
@VV{N_{F'(z)/F(z)}}V  @VV{N_{F'(z)/F(z)}}V  \\
\NNN_{F,z}    @>{d_F}>> (\WWW_M/W_M)^{G_{F(z)}}
\end{CD} \]
and $d_F$ is unique up to $\Hom_{\co[G_K]}(\NNN_{F,z},\WWW_M)$. 
\end{Proposition}
\begin{proof}
See Rubin~\cite{Rubin2}, Proposition 4.4.8, page $87$.
\end{proof}
From Proposition~\ref{proposition:df}, for any element 
$c\in H^0(G_z,\NNN_{F,z}/M\NNN_{F,z})$, let $\tilde c$ be a lifting
of $c$ in $\NNN_{F,z}$.
Then $d_F(N_{F(\mathbf 1)/F})\tilde c$
is an element in $(\WWW_M/W_M)^{G_F}$ and 
$\delta_F d_F(N_{F(\mathbf 1)/F})\tilde c$ is a well defined element in
$H^1(F,W_M)$, independent of the choice of $\tilde c$. 
We denote by
$\kappa$  the map $c\mapsto \delta_F d_F(N_{F(\mathbf 1)/F}) \tilde c$.
Let $\hat d_F$ be a lifting of $d_F$ in $\WWW_M$, then one can see
immediately $\kappa(c)(\gamma)\in W_M$ is exactly $(\gamma-1)
\hat \delta_F d_F(N_{F(\mathbf 1)/F} \tilde c)$.

\subsection{The Kolyvagin recursion}
Fix a number field $F$ of $K$ inside $K_{\infty}$. Fix a infinite 
$\bbz\in Z$ with $x\in X_{F.M}$ for $x\mid \bbz$. 
For every $x\mid \bbz$, one see that
\[ p(x;t)-p(x;1)=Q_x(t)(t-1), \]
for some polynomial $Q_x(t)\in \co[t]$, thus 
$p(x;t)\equiv Q_x(t)(t-1)\pmod M$, 
and $Q_x(t)$ is uniquely determined by the congruence relation since
$t-1$ is not a zero divisor in $\co/M\co[t]$. Fix $Q_x(t)$.

We say that a family 
of classes  
\[ \{\kappa_y\in H^1(F,W_M): y\ \text{finite}, y\mid \bar \bbz\} \]
satisfies the \emph{Kolyvagin recursion} if for every finite $y\mid \bar \bbz$ 
and $x\mid y$, the following formula
\[ Q_x(\Fr^{-1}_x)\kappa_{y/x}(\Fr_x)=\kappa_y
(\sigma_{\bbz(x)})\in W_M \]
holds. We see that Theorem 4.5.4 in Rubin~\cite{Rubin2} essentially
showed that the family $\{\kappa(D_y[\bbz(y)])\}$ satisfies the Kolyvagin
recursion.

\subsection{The universal Kolyvagin recursion implies the Kolyvagin recursion}
We first gather a few
lemmas from Rubin~\cite{Rubin2}:
\begin{Lemma} \label{lemma:zero}
Let $x\in X_{F,M}$, then $p(x;\Fr^{-1}_x)$ annihilates $W_M$.
\end{Lemma}
\begin{proof}
See Rubin~\cite{Rubin2}, Lemma 4.1.2(iv), page $77$.
\end{proof}

\begin{Lemma} \label{lemma:comm}
Let $\hat d_F([z])\in \WWW_M$ be a lifting of $d_F([z])\in \WWW_M/W_M$. 
Then for any $x\mid \bbz$, $\omega$ a prime in $\bar K$ above $x$, 
$g$, $g'$ elements in the decomposition group $\mathcal D$ of $\omega$,
and $\gamma\in G_K$, then
\[ g g'\gamma \hat d_F([z])=g'g\gamma \hat d_F([z]). \]
\end{Lemma}
\begin{proof}
See Rubin~\cite{Rubin2}, Lemma 4.7.1, page $98$.
\end{proof}

\begin{Lemma} \label{lemma:hatdf}
Let  $\hat d_F$ be a lifting
of $d_F$ in $\WWW_M$, then for every $\gamma\in G_K$ and 
$z\mid_s\bbz$ and $x\mid z$, 
\[ N_{G_{\bbz(x)}}\gamma \hat d_F([z])=
p(x;\Fr^{-1}_x)\gamma \hat d_F([z/\bbz(x)]). \]
\end{Lemma}
\begin{proof}
See Rubin~\cite{Rubin2}, Lemma 4.7.3, page $99$.
\end{proof}
\begin{Lemma} \label{lemma:ixzero}
Let $n_x$ be the number of elements in $\co_K/x$. 
Let $r_x(t)=\frac{p(x;t)-p(x;n_x t)}{|G_{\bbz(x)}|}$. Then 
\[ \xi(I_{x,z})\in H^1(F(z)_{\omega}, W_M) \]
for every prime $\omega$ in $F(z)$ above $x$. 
\end{Lemma}
\begin{proof}
See Rubin~\cite{Rubin2}, Corollary 4.8.1, page $102$.
\end{proof}
With the choice of $r_x(t)$ in Lemma~\ref{lemma:ixzero}, 
one then has $\gamma_{\bbz(x)}=p(x; n_x \Fr^{-1}_x)$ and
\begin{Lemma} \label{lemma:gammax}
The map $\gamma_{\bbz(x)}:\NNN_{\bbz/\bbz(x)}\rightarrow 
\NNN_{\bbz/\bbz(x)}$ is injective, thus the Kolyvagin conditions
are satisfied.
\end{Lemma}
\begin{proof} We only need to show that for any finite $z\mid_s \bbz$,
$\gamma_{\bbz(x)}:\NNN_{z/\bbz(x)}\rightarrow \NNN_{z/\bbz(x)}$ is an
injection. In this case, $\Fr^{-1}_x$ induces a linear transformation from
$\NNN_{z/\bbz(x)}\otimes_{\co}\Phi$ to itself, whose eigenvalues
$\lambda$ has (logarithmic) discrete value $0$ at every places above $x$.
 Since $p(x;t)\in \co[t]$, then except the constant term $1$, other terms
in $p(x;n_x \Fr^{-1}_x)$ have discrete value no less than the discrete value of
$n_x$, which is bigger than $0$. Thus any eigenvalue of $\gamma_{\bbz(x)}$
can't be zero and the linear transformation $\gamma_{\bbz(x)}$ is injective.
\end{proof}

Finally we have
\begin{theo}[Kolyvagin recursion]  \label{theo:koly2}
Let $r_x$ be given as Lemma~\ref{lemma:ixzero}.
Let $\{c_y:y\in Y,\ y\mid \bbz\}$ 
be a family of classes in $H^0(G_{\bbz},\NNN_{\bbz}/M\NNN_{\bbz})$
satisfying the universal Kolyvagin recursion related to $\gamma_{\bbz(x)}$.  
Then the family $\{\kappa(c_y):y\in Y,\ y\mid \bbz\}$ satisfies the 
Kolyvagin recursion,
i.e.,
\[ Q_x(\Fr^{-1}_x)\kappa(c_{y/x})(\Fr_x)=\kappa(c_{y})
(\sigma_{\bbz(x)})\in W_M. \]
\end{theo}
\begin{proof} Let $\hat d$ be a lifting of $d_F$ in $\WWW_M$. By the 
definition of the connecting homomorphism $\delta_{F(z)}$, one has
\[ \begin{split}
& \kappa(c_{y/x})(\Fr_x)=(\Fr_x-1)N_{F(\mathbf 1)/F} 
\hat d(c_{y/x})\in W_M,\\
&\kappa(c_{y})(\sigma_{\bbz(x)})=(\sigma_{\bbz(x)}-1)N_{F(\mathbf 1)/F}
\hat d(c_y)\in W_M.
\end{split} \]
Then by Lemmas~\ref{lemma:zero}---~\ref{lemma:gammax}, with the universal
Kolyvagin recursion satisfied by $c_y$ and $c_{y/x}$, 
one has
\[ \begin{split}
 Q_x & (\Fr^{-1}_x)\kappa(c_{y/x})(\Fr_x)-\kappa(c_{y})
(\sigma_{\bbz(x)})\\
&= Q_x(\Fr^{-1}_x)\Fr^{-1}_x\kappa(c_{y/x})(\Fr_x)-\kappa(c_{y})
(\sigma_{\bbz(x)})\\
&=Q_x(\Fr^{-1}_x)(1-\Fr^{-1}_x)N_{F(\mathbf 1)/F}\hat d(c_{y/x})
-(\sigma_{\bbz(x)}-1)N_{F(\mathbf 1)/F}
\hat d(c_y)\\
&=-P(x;\Fr^{-1}_x)N_{F(\mathbf 1)/F}\hat d(c_{y/x})+
\gamma_{\bbz(x)}N_{F(\mathbf 1)/F}\hat d(c_{y/x})\\
&=-|G_{\bbz(x)}|\cdot r_x(\Fr_x^{-1})N_{F(\mathbf 1)/F}\hat d(c_{y/x})\\
&=0,
\end{split} \]
which finishes the proof.
\end{proof}
\begin{rem} 
The proofs of the above Lemmas don't require the use of
the original Kolyvagin classes. Hence we indeed succeed to generalize  
Theorem 4.5.4 in \cite{Rubin2}. We sincerely hope that our more abstract 
construction could lead to pursue new  Euler systems. 
\end{rem}
\begin{rem}
In the special case $T=\ZZ_p(1)$, we see first that
$p(x;t)=1-t$ for every prime $x$. The elements $\kappa(c)\in H^1(F,W_M)$
for $c\in H^0(G,\NNN_{\bbz}/M\NNN_{\bbz})$ are elements inside 
$F^{\times}/F^{\times M}$. Anderson and Ouyang have studied this case
thoroughly in the note ~\cite{AO}.
\end{rem}

\end{document}